\documentclass[12pt]{amsart}
\usepackage{amssymb}
\ifx\mathrm\undefined\let\mathrm\rm\fi
\ifx\mathbf\undefined\let\mathbf\bf\fi
\ifx\mathfrak\undefined\let\mathfrak\frak\fi
\ifx\mathcal\undefined\let\mathcal\cal\fi
\ifx\mathbb\undefined\let\mathbb\Bbb\fi
\ifx\emph\undefined\let\emph\it\fi
\newcommand{\Z}{{\mathbb Z}}

\newcommand{\C}{{\mathbb C}}

\newcommand{\Ref}[1]{{$($\ref{#1}$)$}}
\newcommand{\be}{\begin{displaymath}}
\newcommand{\ee}{\end{displaymath}}
\newcommand{\bea}{\begin{eqnarray*}}
\newcommand{\eea}{\end{eqnarray*}}
\newcommand{\tr}{{\mathrm{tr}}}
\newcommand{\End}{{\mathrm{End}}}
\newcommand{\res}{{\mathrm{res}}}
\newcommand{\h}{{\mathfrak{h}}}

\newcommand{\cR}{{\mathcal{R}}}
\newenvironment{prf}{\noindent{\it Proof\/}:}{$\;\square$
\par\medskip}

\newtheorem%
{thm}{Theorem}[section]
\newtheorem%
{proposition}[thm]{Proposition}
\newtheorem%
{lemma}[thm]{Lemma}
\newtheorem%
{lemmadef}[thm]{Lemma-Definition}
\newtheorem%
{corollary}[thm]{Corollary}
\newtheorem%
{conjecture}[thm]{Conjecture}

\title{Elliptic quantum groups and Ruijsenaars models}
\author[G. Felder and A. Varchenko]{Giovanni Felder${}^{*}$ 
\and Alexander Varchenko${}^{**,1}$}
\thanks{${}^1$Supported in part by NSF grant  DMS-9501290}
\begin{document}
\maketitle
\medskip
\centerline{\it ${}^*$D-MATH, ETH-Zentrum,}
\centerline{\it 8092
Z\"urich, Switzerland}
\medskip
\centerline{\it ${}^{**}$Department of Mathematics,
University of North Carolina at Chapel Hill,}
\centerline{\it Chapel Hill, NC 27599-3250, USA}
\medskip
\centerline{March 1997}
\begin{abstract}
We construct symmetric and exterior powers of the vector
representation of the elliptic quantum groups $E_{\tau,\eta}(gl_N)$.
The corresponding transfer matrices give rise to various integrable
difference equations which could be solved in principle by the
nested Bethe ansatz method. In special cases we recover the
Ruijsenaars systems of commuting difference operators.
\end{abstract}
\section{Introduction}
Elliptic quantum groups \cite{F} are the algebraic structure underlying 
quantum integrable models of statistical mechanics involving elliptic
functions, and the $q$-defor\-ma\-tion of conformal field theory
on elliptic curves.
 
The basic object, appearing in the presentation of elliptic quantum
groups by quadratic relations, is a ``dynamical $R$-matrix'', a
solution of a modification of the Yang--Baxter equation. This
$R$-matrix depends on the spectral parameter and
on an additional parameter lying in the
Cartan subalgebra of a simple Lie algebra. The usual Yang--Baxter
equation is obtained in the limit when the latter parameter tends
to infinity. However, for elliptic solutions, a limit exists only
in trigonometric degenerations. 

In this paper we start a study of the representation
theory of elliptic quantum groups of type $A_{N-1}$, and give some
applications. We construct, using the fusion procedure
\cite{KRS,C,JKMO}, analogues of the symmetric and exterior
powers of the vector representation of $gl_N$.
 Out of these representations, and their
tensor products, one can construct new dynamical $R$-matrices. Taking
partial traces of these $R$-matrices gives rise to families
of commuting difference operators. In a special case, we
recover the Ruijsenaars system \cite{R} 
of  difference operators, which is 
the difference (or ``relativistic'') analogue of the Calogero--Moser
integrable system of differential operators. In particular,
we obtain a simpler proof of the commutativity of Ruijsenaars
operators, in the case of integer coupling constant.
For more general representations, 
one obtains some ``spin generalizations'' of 
Ruijsenaars operators.

The transfer matrix associated to the top exterior power
turns out to be related to the quantum determinant, a
difference operator which is central in the ``operator algebra''
associated to the elliptic quantum group.

In previous papers \cite{FV4,FV5}, a dynamical version of the
Bethe ansatz was developed, and the algebraic integrability
of the Ruijsenaars model was proved in the $A_1$ case. It is
likely that a dynamical version of the ``nested'' Bethe ansatz
(a recursive Bethe ansatz for $A_{N-1}$) can be used
to study the eigenvalue problem for these families of commuting
operators.

Let us conclude this introduction by making some comments
on related papers. The representation theory of elliptic quantum
groups was described in \cite{FV3} for the $A_1$-case.
In the trigonometric degeneration, the
Ruijsenaars operators become the Macdonald difference operators.
They are related to the representation theory of finite dimensional
quantum groups, as pointed out by Etingof and Kirillov \cite{EK}.
In fact, in this case, one can construct them as transfer matrices
with
dynamical $R$-matrices without spectral parameters, see
\cite{ABB}. In the elliptic case, Ruijsenaars operators can
also be obtained as transfer matrices associated to the
Sklyanin algebra, as shown by Hasegawa \cite{H}. His construction
appears to be related to ours, but is more complicated due to
the complexity of the vertex-IRF transformation relating the two
approaches. Hasegawa also gives a space of theta functions which
is invariant under the action of the Ruijsenaars operators.

\noindent{\bf Acknowledgement.} We are grateful to Ivan
Cherednik for inspiring discussions on the fusion procedure.
\section{The elliptic quantum group associated to $gl_N$}
We review here the definition of the elliptic quantum
group $E=E_{\tau,\gamma/2}(gl_N)$ 
associated to $gl_N$ \cite{F} (or rather of its representations)
and the construction
of commuting transfer matrices associated to representations
of $E$ \cite{FV4}. 

Let $\h$ be the Cartan subalgebra of $gl_N$. It is the Abelian
Lie algebra of diagonal complex $N\times N$ matrices.
We identify $\h$ with its dual space via the nondegenerate
bilinear form $(x,y)=\tr(xy)$, and with $\C^N$ via the
 orthonormal basis $\omega_j
=\mathrm{diag}(0,\dots,0,1,0,\dots,0)$ with a 1 in the $j$th
position ($j=1,\dots, N$). 
 
A finite dimensional diagonalisable $\h$-module is a complex
finite dimensional vector space $W$ with a weight decomposition 
$W=\oplus_{\mu\in\h^*}W[\mu]$, so that $\h$ acts on $W[\mu]$
by $xw=\mu(x)w$, ($x\in \h$, $w\in W[\mu]$).
For example, the vector representation of $gl_N$ is $V=\C^N$
with standard basis $e_1,\dots,e_N$ and
with non-zero weight spaces $V[\omega_j]=\C\, e_j$, $j=1,\dots, N$

Let us fix a point $\tau$ in the upper half plane and a generic
complex number $\gamma$. Let 
\[
\theta(z)=-\!\sum_{j\in\Z+\frac12}e^{\pi ij^2\tau+2\pi ij(z+\frac12)},
\]
be Jacobi's first theta function.

Let $R(z,\lambda)\in\End(\C^N\otimes\C^N)$ 
be the $R$-matrix of the elliptic quantum
group $E=E_{\tau,\gamma/2}(gl_N)$. It is a function of the spectral
parameter $z\in\C$ and an additional variable
$\lambda=(\lambda_1,\dots,\lambda_N)
\in\h^*$ in the dual of the Cartan subalgebra
of $gl_N$. It is a solution of the dynamical Yang--Baxter equation
\begin{eqnarray*}
\lefteqn{R(z_1-z_2,\lambda-\gamma h^{(3)})^{(12)}
R(z_1-z_3,\lambda)^{(13)}
R(z_2-z_3,\lambda-\gamma h^{(1)})^{(23)}}
\\
 & &
=
R(z_2-z_3,\lambda)^{(23)}
R(z_1-z_3,\lambda-\gamma h^{(2)})^{(13)}
R(z_1-z_2,\lambda)^{(12)},
\end{eqnarray*}
and is ``unitary'':
\[
R(z,\lambda)R(-z,\lambda)^{(21)}=\mathrm{Id}_{\C^N\otimes\C^N}
\]
We adopt a standard notation: for instance,
$R(z,\lambda-\gamma h^{(3)})^{(12)}$ acts on 
a tensor
$v_1\otimes v_2\otimes v_3$ as 
$R(z,\lambda-\gamma\mu_3)\otimes{\mathrm{Id}}$
if $v_3$ has weight $\mu_3$.
The formula for $R$ is
\[
R(z,\lambda)=\sum_{i=1}^NE_{i,i}\otimes E_{i,i}
+\sum_{i\neq j}\alpha(z,\lambda_i-\lambda_j) E_{i,i}\otimes E_{j,j}
+
\sum_{i\neq j}\beta(z,\lambda_i-\lambda_j) E_{i,j}\otimes E_{j,i}.
\]
The functions $\alpha,\beta$ are  ratios of theta functions:
\[
\alpha(z,\lambda)=\frac{\theta(z)\theta(\lambda+\gamma)}
{\theta(z-\gamma)\theta(\lambda)}\,,
\qquad
\beta(z,\lambda)=-\frac{\theta(z+\lambda)\theta(\gamma)}
{\theta(z-\gamma)\theta(\lambda)}\,,
\]
and $E_{i,j}$ is the matrix such that $E_{i,j}e_k=\delta_{j,k}e_i$.
Note that $R(z,\lambda)$ is invariant under the symmetric group
$S_N$ (the Weyl group of $gl_N$), 
in the sense that for any permutation $\sigma$
\begin{equation}\label{eqsym}
R(z,\sigma\cdot\lambda)=\sigma\otimes\sigma\, R(z,\lambda)\,
\sigma^{-1}\otimes\sigma^{-1},
\end{equation}
where $S_N$ acts linearly on $\h$ and $\C^N$ by permutation
of coordinates.

A representation of the elliptic quantum group $E$ (an $E$-module)
is by definition
a pair $(W,L)$ where $W$ is a (say finite-dimensional,
diagonalisable)  $\h$-module and $L(z,\lambda)$
is a meromorphic function with values in $\End_\h(\C^N\otimes W)$
(the endomorphisms commuting with the action of $\h$),
obeying the relations
\begin{eqnarray*}
\lefteqn{R(z_1-z_2,\lambda-\gamma h^{(3)})^{(12)}
L(z_1,\lambda)^{(13)}
L(z_2,\lambda-\gamma h^{(1)})^{(23)}}
\\
 & &
=
L(z_2,\lambda)^{(23)}
L(z_1,\lambda-\gamma h^{(2)})^{(13)}
R(z_1-z_2,\lambda)^{(12)}.
\end{eqnarray*}
An $E$-submodule of an  $E$-module $(W,L)$ is a pair $(W_1,L_1)$
where $W_1$ is an $\h$-submodule of $W$ such that $\C^n\otimes W_1$
is invariant under the action of all the $L(z,\lambda)$ and
$L_1(z,\lambda)$ is the restriction to this invariant
subspace.  $E$-submodules are $E$-modules.

The basic example of an $E$-module is
$(\C^N,L)$ with $L(z,\lambda)=R(z-w,\lambda)$.
 It is called the vector representation
with evaluation point $w$ and is denoted  by $V(w)$.

Other modules can be obtained by taking tensor
products: if $(W_1,L_1)$ and $(W_2,L_2)$ are $E$-modules,
then also $(W_1\otimes W_2,L)$, with $L(z,\lambda)=
L_1(z,\lambda-\gamma h^{(3)})^{(12)}L_2(z,\lambda)^{(13)}$.

The {\em transfer matrix} associated to an $E$-module $(W,L)$ is
a difference operator acting on the space $F(W[0])$
of meromorphic functions
of $\lambda\in\h^*$ with values in the zero-weight space
 of $W$.
It is defined by the formula
\[
T(z)f(\lambda)=\sum_{\mu}
\tr^{(1)}_{V[\mu]}L(z,\lambda)f(\lambda-\gamma\mu).
\]
The trace is over the (one-dimensional) weight spaces of
$V=\C^N$. More explicitly, let us introduce matrix elements
by $L(z,\lambda)\,e_i\otimes v
=\sum_je_j\otimes L_{ji}(z,\lambda)\,v$. Then 
\[
T(z)f(\lambda)=\sum_{i=1}^NL_{ii}(z,\lambda)f(\lambda-\gamma \omega_i).
\]
It follows from the Yang--Baxter equation that the
transfer matrices commute for different values of the spectral parameters.
For tensor products of vector representations one recovers in this
way the transfer matrices of IRF models. 
In this paper we consider another class of modules.

\section{Symmetric and exterior powers of the vector representation}
\label{Ssep}

For any $n=1,2,\dots$, the symmetric group $S_n$ acts on
$(\C^N)^{\otimes n}
=\C^N\otimes\cdots\otimes\C^N$ by permuting the factors. A tensor
in $(\C^N)^{\otimes n}$ is called symmetric if it is invariant
under the symmetric group. We denote by $S^n(\C^N)$ the space
of symmetric tensors. We denote by $\wedge^n(\C^N)$
the $n$th exterior power of $\C^N$. It is the quotient 
of $(\C^N)^{\otimes n}$
by the subspace $J_n(\C^N)$ spanned by the tensors of the
form $\sigma v-\epsilon(\sigma)v$, $\sigma\in S_n$, $v\in(\C^N)^{n}$.
Here $\epsilon(\sigma)$ is the sign of the permutation $\sigma$.

The $R$ matrix is invertible for generic values of the parameters.
It becomes singular at special values of the spectral parameter.
As is well-known in the non-dynamical case, these singularities
are responsible for the reducibility of the generically irreducible
tensor products of evaluation
representations at special evaluation points.
\begin{lemma}\label{l1} Let $\gamma$ and $\lambda\in\h^*$ be generic.
Then $R(z,\lambda)$ is a nonsingular matrix for all $z\neq\pm\gamma$
(modulo $\Z+\tau\Z$).
\begin{enumerate}
\item[(i)] The image of $R(-\gamma,\lambda)$ is $S^2(\C^N)$
\item[(ii)] The kernel of $R^{{\mathrm{reg}}}(\gamma,\lambda)
\equiv\res_{z=\gamma}R(z,\lambda)$ is
$J_2(\C^N)=S^2(\C^N)$
\end{enumerate}
\end{lemma}

\begin{prf}
We have the ``unitarity property'' $R(z,\lambda)R(-z,\lambda)^{(21)}=1$
which implies that $R(z,\lambda)$ is nonsingular unless $z$ or $-z$
is a pole. This occurs only if $z=\pm\gamma$. The other
claims follow easily from the definition of $R$.
\end{prf}

Let us define operators $W_n(z,\lambda)\in\End((\C^N)^{\otimes n})$
($z\in\C^n$, $\lambda\in \h^*$)
recursively by the conditions:
\begin{eqnarray*}
W_1(z,\lambda)&=&1\\
W_{n+1}(z,\lambda)&=&
R(z_1-z_2,\lambda-\gamma\sum_{j=3}^{n+1}h^{(j)})^{(12)}
\cdots R(z_1-z_n,\lambda-\gamma h^{(n+1)})^{(1n)}
\\
 & &
R(z_1-z_{n+1},\lambda)^{(1\,n\!+\!1)}
\bigl(1\otimes W_n(z_2,\dots,z_{n+1},\lambda-\gamma h^{(1)})\bigr).
\end{eqnarray*}
For instance, $W_2(z,\lambda)=R(z_1-z_2,\lambda)$ and $W_3(z,\lambda)$
is the left-hand side of the dynamical Yang--Baxter equation. It
has therefore two expressions as a product of three $R$-matrices.

\begin{figure}
\setlength{\unitlength}{0.006in}%
\begin{picture}(320,240)(-220,320)
\thicklines
\put(120,320){\line( 0, 1){120}}
\put(120,440){\line( 1, 1){120}}
\put(160,320){\line( 0, 1){ 40}}
\put(160,360){\line( 1, 1){ 40}}
\put(200,400){\line( 1, 1){ 40}}
\put(240,440){\line( 0, 1){ 80}}
\put(240,520){\line(-1, 1){ 40}}
\put(200,320){\line( 1, 1){ 40}}
\put(240,360){\line( 0, 1){ 40}}
\put(240,400){\line(-1, 1){ 40}}
\put(200,440){\line( 0, 1){ 40}}
\put(200,480){\line(-1, 1){ 40}}
\put(160,520){\line( 0, 1){ 40}}
\put(240,320){\line(-1, 1){ 80}}
\put(160,400){\line( 0, 1){ 40}}
\put(160,440){\line(-1, 1){ 40}}
\put(120,480){\line( 0, 1){ 80}}
\put(320,320){\line( 1, 1){120}}
\put(440,440){\line( 0, 1){120}}
\put(400,320){\line( 0, 1){ 40}}
\put(400,360){\line(-1, 1){ 40}}
\put(360,400){\line( 0, 1){ 40}}
\put(360,440){\line( 1, 1){ 40}}
\put(400,480){\line( 0, 1){ 40}}
\put(400,520){\line(-1, 1){ 40}}
\put(360,320){\line(-1, 1){ 40}}
\put(320,360){\line( 0, 1){120}}
\put(320,480){\line( 1, 1){ 80}}
\put(440,320){\line( 0, 1){ 80}}
\put(440,400){\line(-1, 1){120}}
\put(320,520){\line( 0, 1){ 40}}
\end{picture}
\caption{Two expressions for $W_4(z,\lambda)$}
\label{fig1}\end{figure}
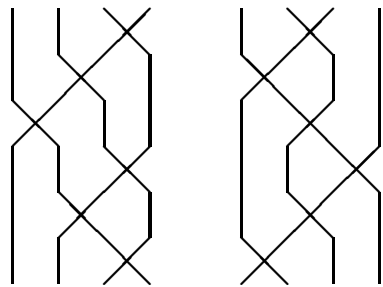

Similarly, $W_n$ can be written in several different ways as
a product of $n(n-1)/2$ $R$-matrices. This is most easily seen
using a graphical representation. The left part of Fig.\  \ref{fig1}
represents the expression  obtained by using the recursive 
definition. 
Let us number the lines
of the drawing from
1 to $n$ from left to right
at the bottom of the drawing. Every crossing between a
line $j$ and a line $k$ represents an $R$-matrix 
\[
R(z_j-z_k,\lambda-\gamma \sum h^{(l)})^{(jk)},
\]
where we adopt the convention that the line numbered $j$ is
at the left of the line $k$ below the crossing. The sum is over
the numbers $l$ assigned to the lines at the left of the
crossing.  The expression
represented by the drawing is then the product of these 
$R$-matrices, the ordering being determined by reading the drawing
from the bottom to the top, so that the bottom crossing corresponds
to the rightmost $R$-matrix.
Using this rule, one can assign a product of $R$-matrices to
any drawing obtained by putting vertical segments and crossings
on top of each other, as in Fig.\  \ref{fig1}, so as to obtain
$n$ lines connecting $n$ points at the bottom with $n$ points
at the top.
A drawing of this kind will be called a diagram. To each
diagram $D$ with $n$ lines
we thus associate a function $W_D(z,\lambda)$ on 
$\C^n\times\h^*$ with values in the endomorphisms of 
$(\C^N)^{\otimes n}$ Every diagram
with $n$ lines induces a permutation of $n$ letters: it sends
$j$ to $k$ if the $j$th point at the bottom is connected to the 
$k$th point at the top.

\begin{lemma}\label{l1a}
If two diagrams $D$, $D'$ 
induce the same permutation then $W_D(z,\lambda)=
W_{D'}(z,\lambda)$.
\end{lemma}
\begin{prf}
We use the fact that the symmetric group $S_n$ is generated by
adjacent transpositions
$s_1,\dots,s_{n-1}$ with relations (a) $s_j^2=1$, (b) $s_js_{j+1}s_j=
s_{j+1}s_js_{j+1}$, (c) $s_js_k=s_ks_j$ if $|j-k|\geq 2$. To a diagram
we associate in an obvious way a word in the generators $s_j$,
so that its image in $S_n$ is the induced permutation.
For example, we associate the word
$s_3s_2s_1s_3s_2s_3$ to the diagram on the left in Fig.\  \ref{fig1}.
If two diagrams induce the same permutation, the corresponding words
can be obtained from each other by applying
a sequence of relations. To each
relation there corresponds a property of $R$-matrices that implies
that the corresponding endomorphisms $W_D(z,\lambda)$ coincide:
namely, we have (a) the ``unitarity'', (b) the dynamical Yang--Baxter
 equation, and (c) $R^{(j,k)}(z,\lambda-\sum_{l\in L} h^{(l)})$
commutes with $R^{(r,s)}(z,\lambda-\sum_{l\in K} h^{(l)})$ if
the sets $\{j,k\}$ and  $\{r,s\}$ are disjoint, and are either
contained in or have empty intersection with $L$ and $K$.
\end{prf}

\begin{corollary}\label{cor}
Let us call a diagram admissible if every line crosses
every other line precisely once. The products of $R$-matrices
associated to any admissible diagram are all equal 
to $W_n(z,\lambda)$.
\end{corollary}

\noindent{\bf Definition.} Let $z^S=(0,\,\gamma,\,2\gamma,\dots,\,
(n-1)\gamma)$, $z^\wedge=((n-1)\gamma\,,\dots,\,2\gamma,\,\gamma,\, 0)$.
 We set
\begin{eqnarray*}
W^S_n(\lambda)&=&W_n(z^S,\lambda),\\
W^\wedge_n(\lambda)&=&\lim_{z\to z^\wedge}\prod_{j=1}^{N-1}
(z_j-z_{j+1}-\gamma)\,W_n(z,\lambda)^{(n,\dots,1)}.
\end{eqnarray*}
\medskip

\noindent{\it Remark.} Note that the spectral parameter in any of the
$R$-matrices of the product defining $W^S_n(\lambda)$
 is always a negative multiple
of $\gamma$, so that there are no divergent $R$-matrices
in this product. Similarly, $W^\wedge_n(\lambda)$ is a product
 of $R$-matrices with spectral parameter $k\gamma$, with $k=2,3,\dots$,
which are finite, and ``regularized'' $R$-matrices
$R^{{\mathrm{reg}}}(\gamma,\lambda)=\res_{z=\gamma}R(z,\lambda)$.

\begin{proposition}\label{p1} Let $\lambda$ be generic.
Then
\begin{enumerate}
\item[(i)] The image of $W^S_n(\lambda)$ is equal to $S^n(\C^N)$.
\item[(ii)] The kernel of $W^\wedge_n(\lambda)$ is equal to $J_n(\C^N)$.
\end{enumerate}
\end{proposition}

\begin{prf} The proofs of the two parts are similar: one inclusion is a 
consequence of Lemma \ref{l1} and Corollary \ref{cor}.
The other inclusion is shown
by counting dimensions at $\gamma=0$. We use the fact
that the dimension of the image (kernel) of a holomorphic family
of matrices at a generic point is at least (at most) the dimension
at a special point.

Let $P$ denote the flip $u\otimes v\mapsto v\otimes u$
on $\C^N\otimes \C^N$. Since every permutation is a product of adjacent
transpositions, we have
\begin{eqnarray*}
S^n(\C^N)&=&\{v\in (\C^N)^{\otimes n}\,|\, P^{(j,j+1)}v=v,\qquad
j=1,\dots, n-1\},\\
J_n(\C^N)&=&\sum_{j=1}^{N-1}\{P^{(j,j+1)}v+v\,|\, v\in
(\C^N)^{\otimes n}\}.
\end{eqnarray*}

\noindent (i) Let us first show that the image of
$W^S_n(\lambda)$ is contained in $S^n(\C^N)$.
It suffices to show that $P^{(j,j+1)}W^S_n(\lambda)=W^S_n(\lambda)$
for all $j=1,\dots,n-1$. This follows from Lemma \ref{l1} (i), and
and the fact that we can always find an admissible  diagram
such that the highest crossing is the one
between line $j$ and line $j+1$ (take any admissible diagram 
and ``move'' lines $j$, $j+1$ so as to push the crossing to the
top).

Let us now show the converse, by considering the limit 
$\gamma\to 0$. The operator $W^S_n(\lambda)$ is
a product of  $R$-matrices whose spectral parameter is
a negative integer multiple of $\gamma$.
Let us write $R_\gamma(z,\lambda)$ for our $R$-matrix to
make the $\gamma$-dependence apparent. The $R$-matrix at 
$z=-k\gamma$ ($k=1,2,\dots$) is then regular as a function of $\gamma$ at
$\gamma=0$ and its limit
\[
\lim_{\gamma\to 0}R_\gamma(-k\gamma,\lambda)=
\frac k{k+1}{\mathrm{Id}}+\frac1{k+1}P,
\]
acts as the identity on $S^2(\C^N)$.
 Hence, if $\gamma=0$, $W^S_n(\lambda)$ 
acts  on  symmetric tensors as the identity
and thus the image contains all symmetric tensors.
It follows that, for generic $\gamma$,
the dimension of the image is at least the dimension of $S^n(\C^N)$.
But since the image is contained in $S^n(\C^N)$, it must coincide with it.

\noindent (ii) We first show that the kernel of $W^\wedge_n(\lambda)$
contains $J_n(\C^N)$. If a vector is of the form $w=P^{(j,j+1)}v+v$,
then $W^\wedge_n(\lambda)$ vanishes on it, as follows from Lemma 
\ref{l1} using a representation of $W^\wedge_n(\lambda)$
corresponding to an admissible diagram such
that the crossing between line $j$ and line $j+1$ is at the bottom.

Let us now show that, in the generic case, the kernel of 
$W^\wedge_n(\lambda)$ is contained in $J_n(\C^N)$.
Let $A_n(\C^N)=\{v\in(\C^N)^{\otimes n}\,|\,P^{(j,j+1)}v=-v,
\qquad j=1,\dots, n\}$ be the space of antisymmetric tensors.
We have a direct sum decomposition $(\C^N)^{\otimes n}=
J_n(\C^N)\oplus A_n(\C^N)$. Indeed, $A_n(\C^N)$ is the orthogonal
complement of $J_n(\C^N)$ with respect to
 the product of standard inner products.

As $\gamma\to 0$, we have
\begin{eqnarray*}
\lim_{\gamma\to 0}R_\gamma(k\gamma,\lambda)
&=&\frac k{k-1}{\mathrm{Id}}-\frac 1{k-1}P,\quad k=2,3,\dots\\
\lim_{\gamma\to 0}\frac1\gamma R_\gamma^{{\mathrm{reg}}}(\gamma ,\lambda)
&=&{\mathrm{Id}}-P.
\end{eqnarray*}
Thus, for $\gamma=0$, $\gamma^{-n+1}W^\wedge_n(\lambda)$ acts as
a non-zero multiple of the identity on $A_n(\C^N)$. The
dimension of the kernel is therefore at most the dimension
of $J_n(\C^N)$, and the claim follows.
\end{prf}

\begin{thm}\label{t1} Let, for $n=0,1,2,\dots$, $V^{\otimes n}(w)$ denote
the $E$-module $V(w)\otimes V(w+\gamma)\otimes\cdots\otimes
V(w+\gamma(n-1))$. Then
\begin{enumerate}
\item[(i)] The subspace $S^n(\C^N)$ is  an
 $E$-submodule
of  $V^{\otimes n}(w)$.
\item[(ii)] The subspace 
$J_n(\C^N)$ is an $E$-submodule of $V^{\otimes n}(w)$.
\end{enumerate}
\end{thm}

\noindent{\bf Definition.} The $E$-module $S^n(\C^N)$ of (i) is called
the $n$th symmetric power of the vector representation with evaluation
point $w$ and is denoted by $S^nV(w)$. The quotient
$V^{\otimes n}(w)/J_n(\C^N)$ by the submodule $J_n(\C^N)$ of (ii) is called
the $n$th exterior power of the vector representation with evaluation
point $w$ and is denoted by $\wedge^nV(w)$. 
\medskip

\noindent What the theorem means is that the $L$-operator on 
$\C^N\otimes(\C^N)^{\otimes\, n}$
\begin{eqnarray}\label{eqL}
L(z,\lambda)&=&R(z\!-\!w,\lambda\!-\!\gamma\sum_{j=2}^n h^{(j)})^{(01)}
R(z\!-\!w\!-\!\gamma,\lambda\!-\!\gamma\sum_{j=3}^n h^{(j)})^{(02)}\\
 & &\cdots
R(z\!-\!w\!-\!\gamma(n\!-\!1),\lambda)^{(0n)},\notag
\end{eqnarray}
(the factors in the tensor products are numbered from $0$ to $n$)
preserves the subspaces $\C^N\otimes S^{n}(\C^N)$ and $\C^N\otimes
J_n(\C^N)$.

In order to prove this theorem, we introduce another operator
on $\C^N\otimes(\C^N)^{\otimes n}$. It is the $L$-operator
corresponding to the ``opposite coproduct'':
\begin{eqnarray*}
 L'(z,\lambda)&=&
R(z-w-\gamma(n-1),\lambda-\gamma\sum_{j=1}^{n-1}h^{(j)})^{(0n)}\\
 & &\cdots 
R(z-w-\gamma,\lambda-\gamma h^{(1)})^{(02)}
R(z-w,\lambda)^{(01)}.
\end{eqnarray*}
\begin{lemma}\label{l3}
\ 
\begin{enumerate}
\item[(i)]
$
L(z,\lambda)\bigl(1\otimes W^S_n(\lambda-\gamma h^{(0)})\bigr)
=\bigl(1\otimes W^S_n(\lambda)\bigr)
 L'(z,\lambda).
$
\item[(ii)]
$
\bigl(1\otimes W^\wedge_n(\lambda)\bigr)
 L(z,\lambda)
=L'(z,\lambda)\bigl(1\otimes W^\wedge_n(\lambda-\gamma h^{(0)})\bigr).
$
\end{enumerate}
\end{lemma}
\begin{prf} 
(i) The left-hand side is, by definition, $W_{n+1}(Z,\lambda)$,
with \[Z=(z,w,w+\gamma,\dots,w+(n-1)\gamma),\] with the notational
convention that the factors are numbered from $0$ to $n$. The
right-hand side is another representation of $W_{n+1}(Z,\lambda)$
as a product of $R$-matrices. It corresponds to the diagram
on the right in Fig.\ \ref{fig1}.   The proof of (ii) is similar:
the claim follows from the identity between two representations
of $W_{n+1}(Z,\lambda)^{(n,\dots,1)}$ in the limit
\[
Z\to (z,w+(n-1)\gamma,\dots,2\gamma,\gamma).
\]
\end{prf}

Theorem \ref{t1} follows from this lemma and Prop.\ \ref{p1}.

{\noindent \bf Example.} The top exterior power $\wedge^NV(w)$ is
a one-dimensional $E$-module. The matrix elements of the
$L$-operator are 
\[
L_{ij}(z,\lambda)=
\delta_{ij}
\frac{\theta(z-w-\gamma)}
{\theta(z-w-\gamma N)}
\prod_{k:k\neq i}
\frac{\theta(\lambda_i-\lambda_k-\gamma)}
{\theta(\lambda_i-\lambda_k)}\, .
\]
\section{Ruijsenaars operators}
The Ruijsenaars operators \cite{R} are integrable 
difference operators that give
a $q$-de\-form\-ation of the Calogero--Moser integrable differential
operators. They act on functions on $\h^*\simeq\C^N$
Let $\Gamma_j$ be the shift by $-\gamma\omega_j$:
\[
\Gamma_jf(\lambda_1,\dots,\lambda_N)
=f(\lambda_1,\dots,\lambda_j-\gamma,\dots,\lambda_N),
\]
and let
$\ell$ be a non-negative integer (the ``coupling constant'').
The corresponding Ruijsenaars operator is (up to conjugation by
a function, see \cite{H})
\[
M=\sum_{i=1}^N\prod_{j:j\neq i}\frac{\theta(\lambda_i-\lambda_j+\ell\gamma)}
{\theta(\lambda_i-\lambda_j)}
\Gamma_i.
\]
It is a symmetric difference operator: if $\sigma\in S_N$ acts
on functions by $\sigma f(\lambda)=f(\sigma^{-1}\lambda)$,
then $\sigma M=M\sigma$ for all $\sigma\in S_N$.

Let us consider the transfer matrix associated to the symmetric
power $S^{n}V(0)$. The zero weight space of this module is trivial
unless $n$ is a multiple of $N$. If $n= N\ell$ then the zero weight
space is one-dimensional and is spanned by the sum of the tensors
$e_{i_1}\otimes\cdots\otimes e_{i_{ N\ell}}$ over all sequences
$(i_j)$ such that each number between $1$ and $N$ occurs precisely 
$\ell$ times. Let us denote this sum by $e$. 
\begin{thm}
Let us identify $S^{ N\ell}(\C^N)[0]$ with $\C$ using the basis $e$,
and let $T(z)$ be the transfer matrix associated to $S^{ N\ell}V(0)$.
Then
\begin{equation}\label{eqT=M}
T(z)=\frac{\theta(z-\gamma\ell)}{\theta(z-\gamma N\ell)}
\,M
\end{equation}
\end{thm}
The rest of this section contains the proof of this theorem.
\begin{lemma}\label{l4}
Let $\bar e$ be the weight zero tensor in $(\C^N)^{\otimes{ N\ell}}$,
\[
\bar e=e_1\otimes\cdots\otimes e_1\otimes\cdots\otimes e_N\otimes\cdots
\otimes e_N,
\]
where each $e_j$ appears $\ell$ times. Let 
\[
g_{N,\ell}(\lambda)=\prod_{1\leq j<k\leq N}\prod_{s=1}^\ell
\frac{\theta(\lambda_j-\lambda_k+\gamma s)}
{\theta(\lambda_j-\lambda_k-\gamma(s-1))}
\, .
\]
Then there exists a non-zero constant $C_{N,\ell}$ such that
\[
W^S_{N\ell}(\lambda)\,\bar e=  C_{N,\ell}\, g_{N,\ell}(\lambda)\, e.
\]
\end{lemma}
\begin{prf} 
We know from Prop.\ \ref{p1} that the left-hand side must be
proportional to the symmetric tensor $e$. It suffices therefore
to compute the coefficient of $\bar e$ in this expression.
It is easy to see that only diagonal elements in the
$R$-matrices give a non-trivial contribution to this coefficient.
Therefore $g_{N,\ell}$ is just a product of functions 
\[
\alpha(z,\lambda)=\frac{\theta(z)\theta(\lambda+\gamma)}
{\theta(z-\gamma)\theta(\lambda)}\, .
\]
The values of $z$ occurring here are of the form $m\gamma$ for some negative
integer $m$, so that no zeros or divergencies appear. The $z$-dependent
factors contribute to $C_{N,\ell}$.

Let us compute the $\lambda$-dependent part of the
 product of functions $\alpha$, using a representation
of $W^S_{N\ell}$ 
by a diagram. As above, we number the lines from 1 to $ N\ell$
from left to right at the bottom of the diagram. Let us
say that the first $\ell$ lines have weight 1, the next $\ell$ lines
have weight 2 and so on. Let us compute the
contribution to the product of a crossing between a line of
weight $j$ and a line of weight $k\neq j$: suppose that
there are $n_i$ lines of weight $i$ to the left of the crossing.
Then the crossing contributes
$\alpha(z,\lambda_j-\lambda_k-\gamma(n_j-n_k))$ (for some $z$)
to the product.
Crossings between lines with
the same weight give a trivial contribution to the product.

Taking together all crossings between lines of weight $j$ and
$k$, and then taking the product over all $j<k$ yields
\[
g_{N,\ell}(\lambda)=
C_{N,\ell}
\prod_{j<k}\;\prod_{r,s=1}^{\ell-1}
\frac
{\theta(\lambda_j-\lambda_k-\gamma(r-s)+\gamma)}
{\theta(\lambda_j-\lambda_k-\gamma(r-s))}\, ,
\]
which can be simplified to the desired expression.
\end{prf}

We  can now complete the proof of the theorem.
For any scalar function $f(\lambda)$, 
\begin{eqnarray*}
T(z)f(\lambda)\, e&=&\sum_j \tr^{(0)}_{V[\omega_j]}L(z,\lambda)\,
e\,
\Gamma_jf(\lambda)\\
 &=&
\sum_j \tr^{(0)}_{V[\omega_j]}L(z,\lambda)
C_{N,\ell}^{-1}
g_{N,\ell}(\lambda-\omega_j)^{-1}
W^S_{N\ell}(\lambda-\gamma\omega_j)\,\bar e\,
\Gamma_jf(\lambda).
\end{eqnarray*}
Since both sides of \Ref{eqT=M} are symmetric difference
operators ($T(z)$ is symmetric as a consequence of \Ref{eqsym})
with shifts by $-\omega_1,\dots,-\omega_N$, 
it is sufficient to show that the coefficients of $\Gamma_1$ on
both sides coincide.
By Lemma \ref{l3},
\begin{eqnarray*}
\tr^{(0)}_{V[\omega_1]}L(z,\lambda)
W^S_{N\ell}(\lambda-\gamma\omega_1)
\,\bar e
&=&
W^S_{N\ell}(\lambda)
\,\tr^{(0)}_{V[\omega_1]}\hat L(z,\lambda)
\,\bar e
\\
 &=&W^S_{N\ell}(\lambda)
\,\bar e\,
\frac{\theta(z-\gamma\ell)}
{\theta(z-\gamma N\ell)}
\prod_{k=2}^N
\frac{\theta(\lambda_1-\lambda_k)}
{\theta(\lambda_1-\lambda_k-\gamma\ell)}\\
 &=&
C_{N,\ell}\,g_{N,\ell}(\lambda)
\,e \,
\frac{\theta(z-\gamma\ell)}
{\theta(z-\gamma N\ell)}
\prod_{k=2}^N
\frac{\theta(\lambda_1-\lambda_k)}
{\theta(\lambda_1-\lambda_k-\gamma\ell)}\, .
\end{eqnarray*}
The second equality is obtained by using the
explicit expression for the $R$-matrices. The calculation
is simple, since only the diagonal entries give a non-vanishing
contribution to the trace, and the product of $R$-matrices
gives a product of functions $\alpha$ which factorize neatly,
as in the proof of Lemma \ref{l4}.

It follows that the coefficient of $\Gamma_1$ in $T(z)$ is
\[
\frac{\theta(z-\gamma\ell)}
{\theta(z-\gamma N\ell)}
\prod_{k=2}^N
\frac{\theta(\lambda_1-\lambda_k)}
{\theta(\lambda_1-\lambda_k-\gamma\ell)}\, 
\frac
{g_{\ell,N}(\lambda)}
{g_{\ell,N}(\lambda-\gamma\omega_1)}\, 
=\,
\frac{\theta(z-\gamma\ell)}
{\theta(z-\gamma N\ell)}
\prod_{k=2}^N
\frac{\theta(\lambda_1-\lambda_k+\gamma\ell)}
{\theta(\lambda_1-\lambda_k)}\, .
\]
We have used here the identity
\[
\frac
{g_{\ell,N}(\lambda)}
{g_{\ell,N}(\lambda-\gamma\omega_1)}\, 
=
\prod_{k=2}^N
\frac{\theta(\lambda_1-\lambda_k+\gamma\ell)
\theta(\lambda_1-\lambda_k-\gamma\ell)}
{\theta(\lambda_1-\lambda_k)^2}\, .
\]
\section{$R$-matrices, higher Ruijsenaars operators and the determinant}

In the preceding section we have considered the $L$-operator
$L(z,\lambda)$
of $S^nV(w)$, the $n$th symmetric powers of the vector representation, 
and its
transfer matrices. In terms of $R$-matrices (see the Appendix), 
$L(z,\lambda)$ is the $R$-matrix
for the $E$-modules $V(z)$ and $S^nV(w)$. More
generally, we may consider the $R$-matrices for
$\wedge^mV(z)$ and $S^nV(w)$. The corresponding transfer matrices
give then a family of commuting difference operators.

Let $V^{\otimes n}(z)=V(z)\otimes V(z+\gamma)\otimes\cdots\otimes
V(z+(n-1)\gamma)$.  Then we have $R$-matrices
$\cR_{V^{\otimes m}(z),V^{\otimes n}(w)}(\lambda)$ 
for $V^{\otimes m}(z)$ and $V^{\otimes n}(w)$ (if $z$, $w$ are generic)
obeying the dynamical Yang--Baxter equation.
 They are products of fundamental $R$-matrices $R(z,\lambda)$,
(cf.\ the Appendix) associated to a diagram where each of the
first $m$ lines crosses each of the last $n$ lines precisely once.

\begin{proposition}
The $R$-matrix $\cR_{V^{\otimes m}(z),V^{\otimes n}(w)}(\lambda)$ 
is invertible for generic $z,w$ and preserves the subspaces
$A\otimes B$ where $A$ is  $(\C^N)^{\otimes m}$,
$S^m(\C^N)$ or $J_m(\C^N)$ and $B$ is
 $(\C^N)^{\otimes n}$,
$S^n(\C^N)$ or $J_n(\C^N)$.
\end{proposition}
\begin{prf}
The invertibility follows from the fact that the product
defining the $R$-matrix contains fundamental $R$-matrices
whose spectral parameter is of the form $z-w+r\gamma$, with
integer $r$, which are invertible for generic $z-w$.
The other claim is proved by commuting the $R$-matrix
with $W^{S,\wedge}_m(\lambda)\otimes{\mathrm{Id}}$ and
${\mathrm{Id}}\otimes W^{S,\wedge}_n(\lambda)$ using Lemma
\ref{l1a} as in the proof of Lemma \ref{l3}.
The details are left to the reader. \end{prf}

In particular we have invertible $R$-matrices
$\cR_{\wedge^mV(z),S^nV(w)}$,
$\cR_{\wedge^mV(z),\wedge^nV(w)}$ obeying dynamical
Yang--Baxter equations for submodules and quotients.
It follows that, if we set $\Gamma_\mu f(\lambda)=f(\lambda-\gamma\mu)$
(so that $\Gamma_j=\Gamma_{\omega_j}$)
then the transfer matrices
\[
T_m(z)=\sum_\mu
\tr_{\wedge^mV(z)[\mu]}\cR_{\wedge^mV(z),S^{N\ell}V(0)}\Gamma_\mu,\qquad 
m=1,\dots N,
\]
are commuting difference operators on the space of functions
of $\lambda\in\h^*$ with values in $S^{N\ell}V(0)[0]\simeq\C$.

Here is an explicit formula for $T_m(z)$ in terms of the
matrix elements of the $L$-operator of $S^nV(0)$: let
$|J|$ denote the cardinality of a subset $J$ of $\{1,\dots,N\}$.
Then
\begin{eqnarray}\label{eqTm}
\lefteqn{
T_m(z)f(\lambda)=
\sum_{1\leq j_1<\cdots<j_m\leq N}
\sum_{\sigma\in S_m}
\epsilon(\sigma)
L_{j_{\sigma(1)}j_1}(z,\lambda-\gamma(\omega_{j_2}
+\cdots+\omega_{j_m}))}\notag
 \\ & &\cdots L_{j_{\sigma(m)}j_m}
(z-(m-1)\gamma,\lambda)
f(\lambda-\gamma(\omega_{j_1}+\cdots+\omega_{j_m})).
\end{eqnarray}
Let us compare these operators commuting with $M$ with
the higher Ruijsenaars operators.

The Ruijsenaars operator $M$ is part of an algebra of $S_N$-symmetric
commuting
difference operators generated by $M_1=M$, $M_2$,\dots, $M_N$, with
\[
M_m=\sum_{J, |J|=m}
\prod_{j\in J,k\neq J}\frac{\theta(\lambda_j-\lambda_k+\gamma\ell)}
{\theta(\lambda_j-\lambda_k)}\prod_{j\in J}\Gamma_j.
\]
\begin{thm} There exist non-zero
meromorphic scalar functions $g_m(z,\gamma)$ so that
\[
T_m(z)=g_m(z,\gamma)M_m,\qquad m=1,\dots, N.
\]
\end{thm}
\begin{prf}
The exterior powers of the vector representation have non-zero weight
spaces
$\wedge^m(\C^N)[\mu]=\C\,e_{j_1}\wedge\cdots\wedge e_{j_m}$,
of weight $\mu=\sum_i \omega_{j_i}$ ($j_1<\cdots<j_m$). It follows
that both sides of the equation we have to prove are $S_N$-symmetric
difference operators
of the form $\sum_{J,|J|=m}A^m_J(z,\lambda)\prod_{j\in J}
\Gamma_j$.
We also know that $T_m(z)$ commutes with $M$. This
implies that the coefficients $A^m_J$ obey difference equations in
$\lambda$.
The idea it to show that these difference equations uniquely determine
the coefficients up to multiplication by a constant.

The coefficients $A_J^m(z,\lambda)$ of $T_m(z)$
are meromorphic functions. As
functions of $\lambda$ they are periodic in each $\lambda_k$ with
period one and they depend only on the differences $\lambda_j-\lambda_k$.
They can be thus considered  as functions of $N-1$ of the
variables $\lambda_j$.

Let $h(x)=\theta(x+\gamma)/\theta(x)$.
If $j\in J$, the fact that $T_m(z)$ commutes with $M$ implies
 the difference equations 
\[
A^m_J(z,\lambda-\gamma\omega_j)=\prod_{k\not\in J}
\frac{h(\lambda_j-\lambda_k-\gamma)}
{h(\lambda_j-\lambda_k)}A^m_J(z,\lambda),\qquad j\in J.
\]
It is easy to check that the coefficients appearing in $M_m$ also obey
these difference equations.
 If $|J|=N$ or $N-1$, this is sufficient to prove the claim.
Namely, in these cases we have a first order
difference equation for all coefficients and all variables $\lambda_j$
(if $|J|=N-1$ we view $A^m_J$ as a function of $N-1$ variables
$\lambda_j$, $j\in J$). We conclude that, for any $J$ with $|J|\geq N-1$,
 the ratio between the coefficient of $\prod_{j\in J}\Gamma_j$
in $T_m(z)$ and the coefficient of $\prod_{j\in J}\Gamma_j$
in $M_m$
is a meromorphic function $g_{J}(z,\lambda;\gamma)$ which, as
a function of $\lambda_k$ ($1\leq k\leq N$), 
is periodic with period $\gamma$
 and with period $1$. Taking $\gamma$ to be an irrational real number,
we see that $g_J$ must be independent of $\lambda$. Since both
$T_m(z)$ and $M_m$ are $S_N$-symmetric, all $g_J$ with $|J|=m$
must be equal to the same function $g_m$.

We have shown that $T_m(z)=g_m(z,\gamma)M_m$ for $m=N$, $N-1$ and some
functions $g_m$ which are not identically zero since they do not
vanish at $\gamma=0$, see below. In particular, this implies that
$T_m(z)$ commutes with
\[
M_{N-1}M_N^{-1}=\sum_{j=1}^N\prod_{k\neq j}h(\lambda_k-\lambda_j)
\Gamma_j^{-1},
\]
since $M_N=\prod_{j=1}^N\Gamma_j$. It follows that, for all $m$, 
the coefficients
$A_J^m$ also obeys the difference equations with respect 
 to the variables $\lambda_j$, $j\not\in J$:
\[
A^m_J(z,\lambda+\gamma\omega_j)=\prod_{k\in J}
\frac{h(\lambda_k-\lambda_j-\gamma)}
{h(\lambda_k-\lambda_j)}A^m_J(z,\lambda),\qquad j\not\in J,
\]
of which the coefficients appearing in $M_m$ are clearly also
a solution.

Proceeding as in the case $m=N, N-1$, we see that the difference
equations uniquely determine all coefficients up to a factor
independent of $\lambda$.
By the $S_N$ symmetry, all factors (for fixed $m$) coincide.

We have still to prove that the functions $g_m(z,\gamma)$ are not
identically zero. As $\gamma\to 0$, $R(z,\lambda;\gamma)$ tends
to the identity, so $T_m(z)$ tends to 
$\sum_{J, |J|=m} \prod_{j\in J}\Gamma_j$. The same holds for $M_m$.
It follows that $g_m(z,0)=1$.
\end{prf}

More generally, for any $E$-module $W$ with $L$-operator 
$L(z,\lambda)$,
we have a family $T_m(z)$, $z\in\C$, $1\leq m\leq N$
of commuting difference operators on $W[0]$-%
valued functions, given by \Ref{eqTm}. 
In general the dimension of $W[0]$ is not
one, so one gets vector valued generalizations of the
 Ruijsenaars model.

These operators are in the ``operator algebra'' \cite{FV1}
of difference operators on $W$-valued functions of $\lambda$.
By definition, the operator algebra is generated by 
the difference operators $\hat L_{ij}(z)=L_{ij}(z,\lambda)\Gamma_j$,
$1\leq i,j,\leq N$. It is clear from \Ref{eqTm} that the
$T_m(z)$ are polynomials in these difference operators.

\newcommand{\Det}{{\mathrm{Det}}}
The transfer matrix $T_N(z)$ associated to the top
exterior power is related to the {\em quantum determinant}
(cf.\ \cite{FV3}, Section 10):
\[
T_N(z)=\frac
{\phi(\lambda-\gamma h)}
{\phi(\lambda)}
\, \Det(z,\lambda)\prod_{j=1}^N\Gamma_j, \qquad
\phi(\lambda)=\prod_{i<j}\theta(\lambda_i-\lambda_j).
\]
The difference operator $\widehat{\Det}(z)=\Det(z,\lambda)\prod\Gamma_j$
defined by this formula is a central element (for all $z$)
of the operator algebra.
This can be seen by writing the Yang--Baxter equation on
$V(z)\otimes \wedge^nV(w)\otimes W$, using the formula in
the example at the end of Section \ref{Ssep} for the 
$R$-matrix of the first two factors.

\section{Conclusions}
We have constructed some finite-dimensional modules over
the elliptic quantum groups associated to $gl_N$.
They are elliptic deformations of
symmetric and exterior powers of the vector representation.
Modules corresponding to more general Young diagrams
will be considered elsewhere.

Transfer matrices associated to modules over elliptic quantum
groups  are commuting difference operators acting on functions
with values in the zero weight space of the quantum space. We have considered
here the special case where the zero weight space is one-dimensional.

In this case the commuting difference operators turn out to be essentially
the Ruijsenaars operators. The advantage of this reformulation
in terms of transfer matrices lies in the fact that we can apply
the Bethe ansatz method to find eigenvectors. In the $gl_2$ case
this was done in \cite{FV4, FV5}. In  general, an adaptation
of the ``nested'' Bethe ansatz to the dynamical case should
give the result. In particular,  one would get a quantum version of the
spectral varieties of \cite{FV1,FV2}.

\appendix

\section{$R$-matrices and transfer matrices}
We give here some details on $R$-matrices and commuting transfer
matrices. In particular, we explain the interpretation of
$R$-matrices as intertwiners for elliptic quantum groups and
give the construction of commuting transfer matrices from
$R$-matrices. The constructions are standard in quantum integrable
systems based on the (non-dynamical) Yang--Baxter equation.
All claims below follow easily from the definitions and from
the dynamical Yang--Baxter equation. 

\newcommand{\Hom}{{\mathrm{Hom}}}
Let $P_{V,W}\in\Hom(V\otimes W,W\otimes V)$ be the flip 
$v\otimes w\mapsto w\otimes v$.
Let $W_1$, $W_2$ be  $E$-modules with $L$-operators $L_1$ and
$L_2$. A {\it morphism}
(or intertwiner) from  $W_1$ to $W_2$
is a meromorphic function $\phi:\h^*\to\Hom_\h(W_1,W_2)$ such
that ${\mathrm{Id}}\otimes
\phi(\lambda)L_1(z,\lambda)=
L_2(z,\lambda){\mathrm{Id}}\otimes\phi(\lambda-\gamma h^{(1)})$. 
An {\it  $R$-matrix} for $W_1$ and $W_2$ is a meromorphic function
$\cR:\h^*\to\End_\h(W_1\otimes W_2)$ so that 
\[
\cR(\lambda)P_{W_2,W_1}: W_2\otimes W_1
\to W_1\otimes W_2,
\]
is a morphism.

For example $\cR(\lambda)=R(z_1-z_2,\lambda)$ is an $R$-matrix
for the evaluation modules $V(z_1)$ and $V(z_2)$, as a consequence
of the dynamical Yang--Baxter equation. Let us call
it the fundamental $R$-matrix with spectral parameter $z_1-z_2$.

Suppose that $W_1$, $W_2$ and $W_3$ are $E$-module and that
$\cR_{W_i,W_j}$ is an  $R$-matrix
for $W_i$ and $W_j$ ($1\leq i<j\leq 3$). Then 
\begin{eqnarray*}
\cR_{W_1\otimes W_2,W_3}(\lambda)
&=&\cR_{W_1,W_3}(\lambda)^{(13)}
\cR_{W_2,W_3}(\lambda-\gamma h^{(1)})^{(23)},
\\
\cR_{W_1,W_2\otimes W_3}(\lambda)
&=&\cR_{W_1,W_3}(\lambda-\gamma h^{(2)})^{(13)}
\cR_{W_1,W_2}(\lambda)^{(12)},
\end{eqnarray*}
are $R$-matrices for $W_1\otimes W_2$, $W_3$ and for 
$W_1$, $W_2\otimes W_3$, respectively.

In particular, we see by iterating this construction that there
are $R$-matrices $\cR_{W_1,W_2}$, obtained as products of fundamental
$R$-matrices, if  $W_1$, $W_2$ are tensor products of vector representations
with generic evaluation points. 
We have, for instance, $\cR_{V(z),W}(\lambda)=L_W(z,\lambda)$, the
$L$ operator of the tensor product $W$ of vector representations.
By construction, these
$R$-matrices obey the dynamical Yang--Baxter equation
\begin{eqnarray*}
\lefteqn{\cR_{W_1,W_2}(\lambda-\gamma h^{(3)})^{(12)}
\cR_{W_1,W_3}(\lambda)^{(13)}
\cR_{W_2,W_3}(\lambda-\gamma h^{(1)})^{(23)}
}
\\ & &= 
\cR_{W_2,W_3}(\lambda)^{(23)}
\cR_{W_1,W_3}(\lambda-\gamma h^{(2)})^{(13)}
\cR_{W_1,W_2}(\lambda)^{(12)},
\end{eqnarray*}
and the ``unitarity'' property
\[
\cR_{W_1,W_2}(\lambda)^{(12)}\cR_{W_2,W_1}(\lambda)^{(21)}=
{\mathrm{Id}}_{W_1\otimes W_2}.
\]
In particular $R$-matrices for generic tensor
products of vector representations are invertible.

Suppose that $\cR_{W_1,W_2}$ is an $R$-matrix for
 $E$-modules $W_1$, $W_2$. The {\it transfer matrix}
with ``auxiliary space'' $W_1$ and ``quantum space'' $W_2$ is the
difference operator acting on functions on $\h^*$ with
values in the zero weight space $W_2[0]$ of $W_2$ 
\[
T_{W_1,W_2} f(\lambda)=\sum_{\mu\in\h^*}\tr_{W_1[\mu]}
(\cR_{W_1,W_2}(\lambda)) \,f(\lambda-\gamma\mu).
\]

\begin{proposition} Suppose that $W_1$, $W_2$, $W_3$
are $E$-modules with $R$-matrices $\cR_{W_i,W_j}$
obeying the dynamical Yang--Baxter equation. Let
for $i=1,2$ $T_i=T_{W_i,W_3}$ be the transfer matrices
with quantum space $W_3$ and assume that $\cR_{W_1,W_2}(\lambda)$
is invertible for generic $\lambda$. Then 
$T_1T_2=T_2T_1$ on $W_3[0]$.
\end{proposition}

\begin{prf} On $W_1\otimes W_2\otimes( W_3[0])$ we can
write the Yang--Baxter equation in the form
\begin{eqnarray*}
\lefteqn{\cR_{W_1,W_3}(\lambda)^{(13)}
\cR_{W_2,W_3}(\lambda-\gamma h^{(1)})^{(23)}
}
\\  & &=
\left(\cR_{W_1,W_2}(\lambda)^{(12)}\right)^{-1}
\cR_{W_2,W_3}(\lambda)^{(23)}
\cR_{W_1,W_3}(\lambda-\gamma h^{(2)})^{(13)}
\cR_{W_1,W_2}(\lambda)^{(12)}.
\end{eqnarray*}
Taking the trace over $W_1\otimes W_2$ yields the result.
\end{prf}

\end{document}